\documentclass[runningheads]{lncse}
\usepackage{multicol} 
\usepackage{makeidx}  

\def\be{\begin{equation}}
\def\ee{\end{equation}}
\def\ba{\begin{eqnarray}}
\def\ea{\end{eqnarray}}
\def\ban{\begin{eqnarray*}}
\def\ean{\end{eqnarray*}}

\def\square#1{\mathop{\mkern0.5\thinmuskip\vbox{\hrule
    \hbox{\vrule\hskip#1\vrule height#1 width 0pt\vrule}\hrule}
    \mkern0.5\thinmuskip}}
\def\Square{\mathchoice{\square{8pt}}{\square{7pt}}{\square{6pt}}{\square{5pt}}}

\newcommand{\EBF}[1]{#1&=&}
\newcommand{\F}[2]{\frac{#1}{#2}}
\newcommand{\LC}{\left\{}
\newcommand{\RC}{\right\}}
\newcommand{\LR}{\left(}
\newcommand{\RR}{\right)}
\newcommand{\LD}{\left.}
\newcommand{\RD}{\right.}
\newcommand{\EC}{\\ & &}
\newcommand{\D}{\partial}
\def\UnitOperator{{\mathchoice {\rm 1\mskip-4mu l} {\rm 1\mskip-4mu l}
{\rm 1\mskip-4.5mu l} {\rm 1\mskip-5mu l}}}

\makeindex
\begin{document}
\mainmatter              
\title{Heat Invariant $E_2$ for Nonminimal
       Operator on Manifolds with Torsion}
\titlerunning{Heat Invariant $E_2$}
%
\author{Vladimir V. Kornyak}
\authorrunning{Vladimir V. Kornyak}  
%
%
\institute{
       Laboratory of Computing Techniques and Automation\\
       Joint Institute for Nuclear Research\\
       141980 Dubna, Russia}

\maketitle
\index{Kornyak@Vladimir V. Kornyak}
\begin{abstract}
Computer algebra methods are applied to investigation of spectral asymptotics
of elliptic differential operators on curved manifolds with torsion and in
the presence of a gauge field.
In this paper we present complete expressions for the second
coefficient $(E_2)$ in the heat kernel expansion for \emph{nonminimal}
operator on manifolds with nonzero torsion.
The expressions were computed for general case of manifolds of arbitrary
dimension $n$ and also for the most important for $E_2$ case $n=2$. The
calculations have been carried out on PC with the help of a program
written in C.
\end{abstract}
\section{Introduction}
Determination of the internal structure of an object via the
spectra of different radiations and waves around the object
is one of the archetypal problems in physics.
More restricted mathematical version of this problem may be
formulated as follows. A manifold (bundle) equipped with such
structures as metric, curvature, torsion, gauge fields etc.
and an elliptic (pseudo)differential operator acting on this
manifold are given. What information about the manifold can one
obtain studying the spectral properties of the operator?
M. Kac phrases the problem in the evocative title of his
paper ``Can one hear the shape of a drum?''\cite{KacM66}.
The answer to this radical question is negative.
In 1964, J. Milnor \cite{Milnor64} found a pair of isospectral
(with respect to the Laplace operator) but non-isometric tori
in dimension sixteen, and recently, in 1999, D. Sch\"uth
\cite{Schueth99} constructed continuous isospectral families
of metrics on the product of spheres
$S^4\times S^3\times S^3.$\footnote{After the Milnor's result many
examples of \emph{multiply connected} isospectral manifolds have
been constructed, but the Sch\"uth's construction is the first example
of closed \emph{simply connected} isospectral but non-isometric
Riemannian manifolds.}
Nevertheless, many global geometric invariants, such as dimension,
volume, and total scalar curvature, are known to be spectrally
determined.
Moreover, various manifolds such as round spheres of dimension
less than or equal to six and 2-dimensional flat tori are uniquely
determined by the spectra of the Laplacian acting on them.

One of the most constructive approaches to study the spectral
properties of operators on manifolds is investigation of
the \emph{heat kernel expansion.} This approach can be described
briefly  as follows. Starting with an elliptic operator $A$
of the order $2r$, acting on a bundle whose base is a compact
close $n$-dimensional manifold $M$,
and introducing an additional ``time'' variable $t$
one can construct the \emph{heat} operator $A - \F{\D}{\D t}$.
Then one can compute the short-time asymptotic expansion of the diagonal
elements of the kernel of this heat operator:
\be
\langle x\vert e^{-tA}\vert x\rangle \sim\sum_{m\geq 0}
\rm E_m(x\vert A) t^{\textstyle \F{m-n}{2r}}, \qquad t\to +0\ .
\label{expansion}
\ee
The coefficients $E_m(x\vert A)$ in this expansion are spectral
invariants of the operator $A$, and encode information about
the asymptotic properties of the spectrum.
These coefficients are called the \emph{heat invariants} or
\emph{heat kernel coefficients}.
They are also widely known under the names Hadamard
coefficients \cite{Hadamard23}\footnote{It was Hadamard
who introduced these coefficients for scalar operator $A$
already in 1923 and established their essential properties.},
Hadamard-Minakshisundaram-DeWitt-Seeley, DeWitt-Seeley-Gilkey
(HMDS or HAMIDEW, DWSG) coefficients, according to papers of these
authors \cite{MinakshisundaramPleijel49,DeWitt65,Seeley67,Gilkey75}.
The heat invariants $E_m$ are of fundamental importance in quantum
field theory, quantum gravity, spectral geometry and topology of manifolds.
Many quantities of interest (such as the effective action,
Green function, anomalies in quantum field theory
\cite{DeWitt65,Seeley67,Gilkey75,BirrellDavies82,BarvinskyVilkovisky85},
the indices of elliptic operators and the invariants of manifolds in
spectral geometry
\cite{AtiyahBottPatodi73,MinakshisundaramPleijel49,KacM66,McKeanSinger67})
are expressed in terms of the heat invariants.

Most papers devoted to computation of the heat invariants deal
with so-called \emph{minimal} operators whose leading term is
a power of the Laplacian and symbol is a scalar (w.r.t space-time indices).
A typical example of such operators is
$$A = -\Square + X.$$
Here $\Square=g^{\mu\nu}D_\mu D_\nu$, $D_\mu$ is a covariant derivative
including generally different connections (affine and spinor connections,
gauge fields), $X$ is a matrix in internal space, i.e., an operator acting
in sections of the bundle. For minimal
operators there are efficient enough methods for computing
the heat invariants.

In this paper we consider \emph{nonminimal} operator of the form
\be
A^{\mu\nu} = -g^{\mu\nu}\Square + a D^{\mu}D^{\nu} + X^{\mu\nu},
\label{nmini}
\ee
where $X^{\mu\nu}$ is a tensor field (bundle indices are assumed implicitly),
$a$ is a scalar parameter
which should satisfy to the condition $a<1$ for the positive
definiteness and, hence, for the ellipticity of operator
(\ref{nmini}).
In recent years special cases of operator (\ref{nmini}) have been
encountered by physicists studying the quantization of gauge and
gravitational fields in arbitrary gauges
\cite{BarvinskyVilkovisky85,BarvinskyVilkovisky87}.
For example, the quantization of Yang-Mills field in an arbitrary
covariant background gauge leads to the operator
$$A^{ab}_{\mu\nu}= -\delta_{\mu\nu}\Square^{ab} -
 \LR\F{1}{\alpha}-1\RR D^{ac}_{\mu}
D^{cb}_{\nu} - 2f^{acb}G^c_{\mu\nu},$$
where  $D_{\mu}$
is a covariant derivative containing the external field potential
$A_{\mu}$, $G_{\mu\nu}$ is a corresponding field strength,
$f^{abc}$ are the structure constants of a corresponding Lie algebra
and $\alpha$ is a scalar (gauge) parameter.
Another example: the quantization of electro-magnetic
field in an external gravitational field leads to the operator
\cite{BarthChristensen83,ChoKantowski95}
$$A_{\mu\nu}=-g_{\mu\nu}\Square - \LR\F{1}{\alpha}
-1\RR D_\mu D_\nu + R_{\mu\nu}$$
(for an analogous operator in quantum gravity see \cite{BarthChristensen83}).

The torsion is defined as antisymmetric part of affine (or linear) connection
$$T^\lambda{}_{\mu\nu} = \Gamma^\lambda{}_{\nu\mu}-\Gamma^\lambda{}_{\mu\nu},$$
where $\Gamma^\lambda{}_{\mu\nu}$ are connection coefficients.
The Einstein's General Relativity (see e.g. \cite{MisnerThorneWheeler73})
is based on a special connection called \emph{Levi-Civita} connection,
i.e., \emph{symmetric}
and compatible with metric affine connection.
This connection can be expressed completely
in terms of metric and is torsionless.
The General Relativity well describes the interaction of the
matter with the gravity as far as macroscopic bulk matter is
considered. However, on the microscopic level, where the elementary
particles posses such quantum property as spin, it seems necessary
to take into account the influence of spin on the geometry of
space-time. To describe the interaction of spinning particles with
the gravitation, a gravitation theory should include the
non-vanishing torsion.
In 1922 Elie Cartan first pointed out \cite{Cartan22,Cartan23-25} that there is
no \emph{a priori} reason to assume an affine connection to be symmetric in the
context of General Relativity. He proposed also a theory of gravitation
with torsion which development is known now as the Einstein-Cartan theory.
The torsion arises naturally in the different (based on Poincar\'e and affine
 groups) gauge theories of gravity developed in the recent years
 (see Refs. \cite{HehlHeydeKerlickNester76,Hehl80}).
Moreover, all kinds of modern superstring theories \cite{GreenSchwarzWitten87}
(for the recent review see e. g. \cite{Kiritsis97}), which allow to deduce
the properties of space-time, also predict, along with the metric, the existence of
torsion.

In \cite{GusyninKornyak99} we computed the heat invariants for
operator (\ref{nmini}) up to $E_4$ but for manifolds without
torsion. In this paper we consider more general (and computationally
much more difficult) case of manifolds with torsion.

\section{Algorithm and Implementation}
The algorithm we use was developed by V. Gusynin
\cite{Gusynin89,Gusynin90}. This algorithm is based on the
covariant generalization of the \emph{pseudodifferential calculus}
given by Widom \cite{Widom80}. The main advantage of this
algorithm is its universality. It can be applied to the wide
class of pseudodifferential operators, in particular, to the
nonminimal and higher-order operators intractable by other
methods such as the \emph{DeWitt ansatz} \cite{DeWitt65}
for heat kernel matrix elements.

The algorithm has the following main features.
For a positive elliptic operator $A$ the spectrum of which lies
inside a contour $C$, the heat operator $\exp(-tA)$ can be
expressed in terms of the resolvent
$(A-\lambda)^{-1}$ via the formula
\be
e^{-tA}=\int_{C}\F{id\lambda}{2\pi}e^{-t\lambda}(A-\lambda)^{-1}.
\ee
The pseudodifferential calculus method uses the following
representation for the matrix elements of the resolvent
\be
G(x,x',\lambda)\equiv
        \langle x\vert\F{1}{A-\lambda}\vert x'\rangle
        =\int\F{d^{n}k}{(2\pi)^{n}\sqrt{g(x')}}e^{il(x,x',k)}
                                        \sigma(x,x',k;\lambda),
\ee
where $\sigma(x,x^\prime,k;\lambda)$ is an amplitude,
$l(x,x^\prime,k)$ is a (real) phase function which is a biscalar
with respect to general coordinate transformations,
$k$ is a wave vector.

The resolvent satisfies the equation $(A-\lambda)G=1$
which leads to the equation for the amplitude:
\be
(A(x,D_{\mu}+iD_{\mu}l)-\lambda)\sigma(x,x',k;\lambda)=I(x,x'),
\label{ampleq}
\ee
where $I(x,x')$ is a transport function having both bundle and Lorentz
indices.

In the pseudodifferential calculus, it is assumed that in the flat space
the phase function has the form $l = (x-x')_\mu k^\mu.$
The covariant analogue of the linearity of the function $l$
is based on the requirement that all higher-order symmetrized covariant
derivatives of $l(x,x',k)$ vanish at the points $x = x'$,
i. e., satisfy the infinite set of relations \cite{Widom80}:
\be
[\{D_{\mu_{1}}\ldots D_{\mu_{m}}\}l]=0,\qquad m > 1,
\label{lcond}
\ee
where $\{\ldots\}$ means symmetrizing in all indices,
and $[\ldots]$ means transition to coincidence limit $(x=x')$.
In an analogous way, the covariant transport function should
satisfy the relations:
\be
 [\{D_{\mu_{1}}\ldots D_{\mu_{m}}\}I]=0,\quad m\geq 1.
\label{Icond}
\ee
Equations (\ref{lcond}) and (\ref{Icond})  together with the ``initial conditions"
$[l]=0, [D_{\mu}l]=k_{\mu}$ and $[I]=\UnitOperator$ (unit operator)
allow one to compute the coincidence limits for  nonsymmetrized
covariant derivatives
$[D_{\mu_{1}} \ldots D_{\mu_{m}} l]$ and
$[D_{\mu_{1}} \ldots D_{\mu_{m}} I]$.
These nonsymmetrized derivatives
are obtained directly from (\ref{lcond}) and (\ref{Icond}) by reducing
all terms to a unified index ordering with the help of the Ricci identity.
The resulting expressions are universal polynomials
in the torsion $T^\lambda{}_{\mu\nu}$, curvature tensor
$R^\lambda{}_{\mu\nu\eta}$, gauge curvature $W_{\mu\nu}$ and
their covariant derivatives.
In fact, once computed and stored the coincidence limits
$[D_{\mu_{1}} \ldots D_{\mu_{m}} l]$ and $[D_{\mu_{1}} \ldots D_{\mu_{m}} I]$
can be used in many calculations for different operators $A$.
The functions $l(x,x^\prime,k)$ and $I(x,x^\prime)$, introduced with
the help of formulas (\ref{lcond}) and (\ref{Icond}),\footnote{The existence
 of these functions has been proved in \cite{Widom80}.} play an important
 role in the covariant pseudodifferential calculus called also
\emph{intrinsic symbolic calculus} \cite{Widom80}.
In fact, just these universal functions
manifest the geometric properties of a base manifold and a bundle.

Expanding the amplitude $\sigma$ in degrees of homogeneity of
$k$: $$\sigma=\sum_{m=1}^{\infty}\sigma_{m}(x,x',k;\lambda),$$
we obtain the recursion equations for
$\sigma_m$ from equation (\ref{ampleq}).
For example, for operator (\ref{nmini}) these recursion expressions take the form
\ban
&&A^{\mu\lambda}\sigma_{0\lambda\nu}=I^\mu_\nu,\\
&&A^{\mu\lambda}\sigma_{1\lambda\nu}+i\left[-g^{\mu\lambda}(\Square
l+2D^\eta lD_\eta)+a(D^\mu D^\lambda l+D^\mu lD^\lambda+D^\lambda lD^\mu) \right]\\
&&\quad \times\ \sigma_{0\lambda\nu}=0,\\
&&\ \vdots \\
&&A^{\mu\lambda}\sigma_{m\lambda\nu}+i\left[-g^{\mu\lambda}(\Square l+2D^\eta lD_\eta)
+a(D^\mu D^\lambda l+D^\mu lD^\lambda+D^\lambda lD^\mu)\right]\\
&&\quad \times\ \sigma_{(m-1)\lambda\nu}+(-g^{\mu\lambda}\Square+aD^\mu D^\lambda
+X^{\mu\lambda})\sigma_{(m-2)\lambda\nu}=0,\quad
m\geq 2,
\ean
where the matrix
\ban
A^{\mu\nu}=g^{\mu\nu}(D^\eta lD_\eta l-\lambda)-aD^\mu lD^\nu l
\ean
is the principal symbol for operator (\ref{nmini}).
Solving the recursion equations we obtain $\sigma_m$.
The heat invariants are expressed in terms of integrals of
the coincidence limits
$[\sigma_m]$:
\be
E_{m}(x\vert A)=\int\F{d^{n}k}{(2\pi)^{n}\sqrt{g}}
  \int_{C}\F{id\lambda}{2\pi}e^{-\lambda}[\sigma_{m}](x,k,\lambda)\equiv
  J([\sigma_m]).
\label{eqforem}
\ee
The integrals in (\ref{eqforem}) can be expressed in terms of
gamma and Gauss hypergeometric functions for a wide class of
operators $A$.
The typical integral of terms of the coincidence limit
$[\sigma_m]$ takes the form
\ban
\lefteqn{J\Bigl(\F{k^{2p}k_{\mu_{1}}\ldots k_{\mu_{2s}}}
 {(k^{2r}-\lambda)^l[(1-a)k^{2r}-\lambda]^{m}}\Bigr)=}\\
 & & g_{\{\mu_{1}\ldots \mu_{2s}\}}\F{\Gamma((p+s+n/2)/r)}
        {(4\pi)^{n/2}2^{s}r\Gamma(n/2+s)\Gamma(l+m)}
        F(m,(p+s+n/2)/r;l+m;a),
\label{integral}
\ean
where $g_{\{\mu_1\ldots \mu_{2s}\}}$ is a symmetrized sum of products
of metric tensors.
Using the fact that $m$ and $l$ are whole numbers, one can
express the hypergeometric function in (\ref{integral})
in terms of elementary functions with the help
of the Gauss relation
\be
a(1-z)F(a+1,b;c;z)=(c-a)F(a-1,b;c;z)+(2a-c-az+bz)F(a,b;c;z),
\label{Gauss}
\ee
and using then the formula \cite{PrudnikovBrychkovMarichev86}
\ban
F(1,b;m;z)=(m-1)!\F{(-z)^{1-m}}{(1-b)_{m-1}}\Bigl[(1-z)^{m-b-1}-
\sum_{k=0}^{m-2}\F{(b-m+1)_k}{k!}z^k\Bigr],\\
m=1,2,\ldots ,\quad m-b\neq 1,2,\ldots,
\label{Prud}
\ean
where $(a)_k=a(a+1)\ldots(a+k-1)$ is the Pochhammer symbol
(shifted factorial).
During simplification of tensor expressions we use various
symmetry properties of the tensors $R^\lambda{}_{\eta\mu\nu},
T^\lambda{}_{\mu\nu}, W_{\mu\nu}$,
and also the Ricci identity
\ban
&&[D_{\mu},D_{\nu}]\varphi^{\eta_{1}\ldots \eta_{l}}_{\lambda_{1}\ldots \lambda_{k}}
=
\sum_{i=1}^{l}R^{\eta_{i}}{}_{\alpha\mu\nu}\varphi^{\eta_{1}\ldots \eta_{i-1}\alpha
\eta_{i+1}\ldots \eta_{l}}_{\lambda_{1}\ldots \lambda_{k}}
\\
&&\quad -\sum_{i=1}^{k}R^{\alpha}{}_{\lambda_{i}\mu\nu}\varphi^{\eta_{1}\ldots
\eta_{l}}_{\lambda_{1}\ldots
\lambda_{i-1}\alpha \lambda_{i+1}\ldots \lambda_{k}}
+\ T^{\alpha}{}_{\mu\nu}D_{\alpha}\varphi^{\eta_{1}\ldots
\eta_{l}}_{\lambda_{1}\ldots \lambda_{k}}
+ W_{\mu\nu}\varphi^{\eta_{1}\ldots \eta_{l}}_{\lambda_{1}\ldots \lambda_{k}},
\ean
the Bianchi identities for both affine and gauge curvatures
\ban
&&
D_{\alpha}R^{\beta}{}_{\gamma\delta\epsilon}
+D_{\delta}R^{\beta}{}_{\gamma\epsilon\alpha}
+D_{\epsilon}R^{\beta}{}_{\gamma\alpha\delta}
\\
&&\qquad\qquad
+T^{\lambda}{}_{\alpha\delta}R^{\beta}{}_{\gamma\epsilon\lambda}
+T^{\lambda}{}_{\delta\epsilon}R^{\beta}{}_{\gamma\alpha\lambda}
+T^{\lambda}{}_{\epsilon\alpha}R^{\beta}{}_{\gamma\delta\lambda}=0,
\\
&&
D_{\alpha}W_{\beta\gamma}+D_{\beta}W_{\gamma\alpha}
+D_{\gamma}W_{\alpha\beta}
\\
&&\qquad\qquad
+W_{\alpha\lambda}T^{\lambda}{}_{\beta\gamma}
+W_{\beta\lambda}T^{\lambda}{}_{\gamma\alpha}
+W_{\gamma\lambda}T^{\lambda}{}_{\alpha\beta}=0,
\ean
and the cyclic identity
\ban
&&R^{\alpha}_{\beta\gamma\delta}+R^{\alpha}_{\gamma\delta\beta}
+R^{\alpha}_{\delta\beta\gamma}
+D_{\beta}T^{\alpha}{}_{\gamma\delta}+D_{\gamma}T^{\alpha}{}_{\delta\beta}
+D_{\delta}T^{\alpha}{}_{\beta\gamma}
\\
&&\qquad\qquad
+T^{\alpha}{}_{\beta\lambda}T^{\lambda}{}_{\gamma\delta}
+T^{\alpha}{}_{\gamma\lambda}T^{\lambda}{}_{\delta\beta}
+T^{\alpha}{}_{\delta\lambda}T^{\lambda}{}_{\beta\gamma}=0.
\ean

The above algorithm has been implemented in the C language.
The C code of total length about 11000 lines contains about 250
functions for different manipulations with tensors and scalars.
These functions are gathered into two programs DWSGCOEF and COLIM.

The COLIM program computes coincidence limits of the $l(x,x',k)$ and
$I(x,x')$
functions and writes them to the disk. Once computed and
stored\footnote{For the operators of different tensor ranks
$A, A^{\mu\nu},\ldots$
the coincidence limits for the functions
$I, I^{\mu\nu},\ldots$ should be computed separately.}
the coincidence limits, being universal functions,
can be used in many calculations for different operators $A$.

The DWSGCOEF program computes $E_m$ coefficients by the following
steps:
\begin{enumerate}
 \item  \emph{Reading input information (operator, order $m$, etc.)}
 \item  \emph{Computing a set of asymptotic operators for constructing
   recursion equations.}
 \item\label{SigmaM}  \emph{Computing $\sigma_m$ with the help of the recursion
   equations.}
 \item\label{TakeColim}  \emph{Taking the coincidence limit $[\sigma_m]$.}
 \item  \emph{Integrating $[\sigma_m]$ to obtain the coefficient $E_m$.}
 \item\label{ReplaceColim}  \emph{Substituting tensor expressions for
   \([D_{\mu_{1}} \ldots D_{\mu_{k}} l]\) and
   \([D_{\mu_{1}} \ldots D_{\mu_{k}} I]\) into $E_m$.}
 \item  \emph{Reducing hypergeometric to elementary functions in the scalar
   coefficients ( $C_i$ in the formulas of Section
   \ref{e2tnminiout})
   including in the heat invariants in the case
   of nonminimal or higher-order operator. Eliminating possible linear
   dependencies   among these scalar coefficients\footnote{Usually
   there are many dependencies among the scalar coefficients which
   are not seen in terms of hypergeometric functions, i.e., quite different
   hypergeometric expressions may be reduced sometimes to the same elementary function.}
   to make the resulting formulas  as compact as possible.}
 \item  \emph{Output $E_m$ (and its Lorentz trace in the nonminimal case).}
\end{enumerate}
To cut down the swelling of the intermediate expressions,
we use \emph{term-by-term\/}
strategy, i.e., the most cumbersome Steps
\ref{TakeColim}-\ref{ReplaceColim} are applied
consecutively to single terms of $\sigma_m$ generated
 during the execution of Step \ref{SigmaM}.
\section{Heat Invariant $E_2$}\label{e2tnminiout}
We present here the full expression for the coefficient
$E_2$ and also its trace with respect to Lorentz indices
for nonminimal operator (\ref{nmini})
on a curved manifold with the torsion and gauge field manifested itself in the
gauge curvature $W_{\mu\nu}$. We consider the case of arbitrary dimension
$n$ and also the most important\footnote{The \emph{Atiyah--Singer index}
of an elliptic operator on a manifold of
dimension $n$ can be expressed in terms of
an integral of $E_n$ over the manifold.} for $E_2$ case $n=2.$
In the below formulas the indices $\alpha, \beta, \gamma$ and
$\mu, \nu$ are dummy and free, correspondingly.
We use the following definition for the torsion trace:
$T_\mu=T^\alpha{}_{\alpha\mu}$.
\subsection{Full Expression}
\ban
\EBF{E_2}
(4\pi)^{-\F{n}{2}}\LC-C_1X^{\mu\nu}
-C_2\LR X^{\nu\mu}+g^{\mu\nu}X_\alpha{}^\alpha\RR
+C_3\LR W^{\mu\nu}+\F{8}{3}D_\alpha T^{\alpha\mu\nu}\RD\RD
\EC
\LD+\F{19}{6}T_\alpha T^{\alpha\mu\nu}\RR
+C_4R^{\mu\nu}-C_5R^{\nu\mu}
+C_6\LR D_\alpha T^{\mu\alpha\nu}+D_\alpha T^{\nu\alpha\mu}\RR
\EC
+C_7T_{\alpha\beta}{}^\mu T^{\alpha\beta\nu}
+C_8T_{\alpha\beta}{}^\mu T^{\beta\alpha\nu}
+C_9\LR T_{\alpha\beta}{}^\mu T^{\nu\alpha\beta}
+T_{\alpha\beta}{}^\nu T^{\mu\alpha\beta}\RR
\EC
-C_{10}T^\mu{}_{\alpha\beta}T^{\nu\alpha\beta}
+C_{11}D^\mu T^\nu
-C_{12}D^\nu T^\mu
+C_{13}T_\alpha\LR T^{\mu\alpha\nu}+T^{\nu\alpha\mu}\RR
\EC
-C_{14}T^\mu T^\nu
+C_{15}g^{\mu\nu}\LR R+D_\alpha T^\alpha\RR
+C_{16}g^{\mu\nu}T_\alpha T^\alpha
\EC
\LD
+\LR C_{16}-C_{15}\RR
g^{\mu\nu}T_{\alpha\beta\gamma}T^{\beta\alpha\gamma}
-C_{17}g^{\mu\nu}T_{\alpha\beta\gamma}T^{\alpha\beta\gamma}
\RC.
\ean
Coefficients $C_i$ in arbitrary dimension $n$:
\ban
C_1&=&\F{1}{a(n-2)n(n+2)}\LC(1\!-\!a)^{-\F{n}{2}}\LR-3an-6a+4n+4\RR\RD+an^3-2an^2\\
& &\LD-3an+6a-4n-4\RC,\\
C_2&=&\F{1}{a(n-2)n(n+2)}\LC(1\!-\!a)^{-\F{n}{2}}(an+2a-4)+an-2a+4\RC,\\
C_3&=&\F{1}{a(n-2)n}\LC(1\!-\!a)^{1-\F{n}{2}}(an-8)+3an-8a+8\RC,\\
C_4&=&\F{1}{6a(n-2)n(n+2)}\LC(1\!-\!a)^{-\F{n}{2}}\LR17a^2n^2+34a^2n-17an^2
-168an\RD\RD\\
& &\LD\LD-268a+140n+256\RR
-53an^2+40an+268a-140n-256\RC,\\
C_5&=&\F{1}{6a(n-2)n(n+2)}\LC(1\!-\!a)^{-\F{n}{2}}
\LR15a^2n^2+30a^2n-15an^2-152an\RD\RD\\
& &\LD\LD-244a+116n+256\RR-43an^2+24an+244a-116n-256\RC,\\
C_6&=&\F{1}{6a^2(n-2)n(n+2)}\LC(1\!-\!a)^{-\F{n}{2}}
\!\LR a^3n^2\!+\!2a^3n\!-\!a^2n^2\!-\!20a^2n\!-\!36a^2\RD\RD\\
& &\LD\LD+12an+96a-48\RR+a^2n^2-16a^2n+36a^2+12an-96a+48\RC,\\
C_7&=&\F{1}{6a^2(n-2)n(n+2)(n+4)}\LC(1\!-\!a)^{-\F{n}{2}}
\LR-a^3n^3-6a^3n^2+a^2n^3-\RD\RD\\
& &8a^3n+30a^2n^2+152a^2n\!-\!24an^2+192a^2\!-\!240an\!-\!576a+144n\\
& &\LD+288\RR\!-\!7a^2n^3+18a^2n^2+64a^2n\!-\!48an^2\!-\!192a^2+96an+576a\\
& &\LD-144n-288\RC,\\
C_8&=&\F{1}{6a^2(n-2)n(n+2)}\!\LC(1\!-\!a)^{1-\F{n}{2}}
\LR a^2n^2+2a^2n-24an-48a+144\RR\RD\\
& &\LD-7a^2n^2+34a^2n-48a^2-48an+192a-144\RC,\\
C_9&=&\F{1}{a^2(n-2)n(n+2)(n+4)}\LC(1\!-\!a)^{-\F{n}{2}}
\LR a^2n^2+6a^2n+8a^2-12an\RD\RD\\
& &\LD\LD-48a+48\RR-a^2n^2+6a^2n-8a^2-12an+48a-48\RC,\\
C_{10}&=&\F{1}{12a^2(n-2)n(n+2)(n+4)}\LC(1\!-\!a)^{-\F{n}{2}}
\LR-a^3n^3-6a^3n^2+a^2n^3-\RD\RD\\
& &\LD8a^3n+30a^2n^2+152a^2n+192a^2-192an-768a+576\RR-a^2n^3\\
& &\LD-6a^2n^2+88a^2n-192a^2-96an+768a-576\RC,\\
C_{11}&=&\F{1}{3a^2(n-2)n(n+2)}\LC(1\!-\!a)^{1-\F{n}{2}}
\LR-9a^2n^2-18a^2n+76an+152a\RD\RD\\
& &\LD\LD-24\RR+2\LR-13a^2n^2+6a^2n+76a^2-32an-88a+12\RR\RC,\\
C_{12}&=&\F{1}{3a^2(n-2)n(n+2)}\LC2(1\!-\!a)^{1-\F{n}{2}}
\LR-5a^2n^2-10a^2n+38an+76a\RD\RD\\
& &\LD+12\RR\LD-31a^2n^2+26a^2n+152a^2-88an-128a-24\RC,\\
C_{13}&=&\F{1}{6a(n-2)n(n+2)(n+4)}\LC(1\!-\!a)^{-\F{n}{2}}
\LR a^3n^3+6a^3n^2-a^2n^3+8a^3n\RD\RD\\
& &\LD-18a^2n^2\!-\!80a^2n\!+\!12an^2\!-\!96a^2\!+\!72an\!+\!96a\!-\!48n\!+\!96\RR\!+\!a^2n^3\\
& &\LD-18a^2n^2+8a^2n+12an^2+96a^2-120an-96a+48n-96\RC,\\
C_{14}&=&\F{1}{2a^2(n-2)n(n+2)}\LC(1\!-\!a)^{1-\F{n}{2}}
\LR-a^2n^2-2a^2n+8an+16a-16\RR\RD\\
& &\LD-a^2n^2-2a^2n+16a^2-32a+16\RC,\\
C_{15}&=&\F{1}{6a(n-2)n(n+2)}\LC(1\!-\!a)^{-\F{n}{2}}
\LR\!-\!a^2n^2\!-\!2a^2n\!+\!an^2\!+\!8an\!+\!12a\!-\!24\RR\RD\\
& &\LD+an^3-an^2-12a+24\RC,\\
C_{16}&=&\F{1}{12a^2(n-2)n(n+2)(n+4)}\LC(1\!-\!a)^{-\F{n}{2}}
\LR-a^3n^3-6a^3n^2+a^2n^3\RD\RD\\
& &\LD-8a^3n+6a^2n^2+8a^2n+48an+192a-288\RR+a^2n^4+3a^2n^3\\
& &\LD+2a^2n^2-48a^2n+96an-192a+288\RC,\\
C_{17}&=&\F{1}{24a^2(n-2)n(n+2)(n+4)}\LC(1\!-\!a)^{-\F{n}{2}}
\LR-a^3n^3-6a^3n^2+a^2n^3\RD\RD\\
& &\LD-8a^3n+30a^2n^2\!+\!152a^2n+192a^2\!-\!192an\!-\!768a\!+\!576\RR\!+\!a^2n^4\\
& &\LD+3a^2n^3-10a^2n^2+72a^2n-192a^2-96an+768a-576\RC.
\ean
Coefficients $C_i$ in the dimension $n=2$:
\ban
C_1&=&-\F{3\ln(1\!-\!a)}{4a}-\F{2-a}{8(1\!-\!a)},\\
C_2&=&\F{\ln(1\!-\!a)}{a}+\F{2-a}{2(1\!-\!a)},\\
C_3&=&-\F{(a-4)\ln(1\!-\!a)}{2a}+2,\\
C_4&=&\F{(17a-67)\ln(1\!-\!a)}{12a}+\F{137a-134}{24(1\!-\!a)},\\
C_5&=&\F{(15a-61)\ln(1\!-\!a)}{12a}+\F{119a-122}{24(1\!-\!a)},\\
C_6&=&\F{(a^2-9a+6)\ln(1\!-\!a)}{12a^2}+\F{(3a-2)(a-2)}{8a(1\!-\!a)},\\
C_7&=&-\F{(a^2-12a+12)\ln(1\!-\!a)}{12a^2}-\F{2a^2-9a+6}{6a(1\!-\!a)},\\
C_8&=&-\F{(a^2-12a+18)\ln(1\!-\!a)}{12a^2}+\F{a-6}{4a},\\
C_9&=&\F{(a-2)\ln(1\!-\!a)}{4a^2}-\F{a^2-12a+12}{24a(1\!-\!a)},\\
C_{10}&=&-\F{(a^2-12a+12)\ln(1\!-\!a)}{24a^2}-\F{2a^2-9a+6}{12a(1\!-\!a)},\\
C_{11}&=&\F{(9a^2-38a+3)\ln(1\!-\!a)}{6a^2}-\F{73a-6}{12a},\\
C_{12}&=&\F{(10a^2-38a-3)\ln(1\!-\!a)}{6a^2}-\F{79a+6}{12a},\\
C_{13}&=&\F{(a-6)\ln(1\!-\!a)}{12a}+\F{2a-3}{6(1\!-\!a)},\\
C_{14}&=&\F{(a^2-4a+2)\ln(1\!-\!a)}{4a^2}-\F{3a-2}{4a},\\
C_{15}&=&-\F{(a-3)\ln(1\!-\!a)}{12a}-\F{7a-10}{24(1\!-\!a)},\\
C_{16}&=&-\F{(a^2-6)\ln(1\!-\!a)}{24a^2}-\F{3a^2+a-6}{24a(1\!-\!a)},\\
C_{17}&=&-\F{(a^2-12a+12)\ln(1\!-\!a)}{48a^2}-\F{3a^2-10a+6}{24a(1\!-\!a)}.
\ean
\subsection{Lorentzian Trace}
\ban
\EBF{{\rm tr_L} E_2}
(4 \pi)^{-\F{n}{2}}\LC-C_1X_\alpha{}^\alpha
+ C_2R-C_3T_{\alpha\beta\gamma}T^{\alpha\beta\gamma}
+C_4T_{\alpha\beta\gamma}T^{\beta\alpha\gamma}\RD\\
& &
\LD+C_5D_\alpha T^\alpha+\LR C_4+C_5\RR T_\alpha T^\alpha\RC
\ean
$C_i$ for the trace in arbitrary dimension $n$:
\ban
C_1&=&\F{(1\!-\!a)^{-\F{n}{2}}+n-1}{n},\quad C_2\ =\ \F{(1\!-\!a)^{-\F{n}{2}}
(-an+n+6)+n^2-n-6}{6n},\\
C_3&=&\F{1}{24an(n+2)}\LC-(1\!-\!a)^{-\F{n}{2}}\LR a^2n^2+2a^2n-an^2-26an-48a+96\RR\RD\\
& &\LD+an^3+an^2+22an-48a+96\RC,\\
C_4&=&\F{1}{12an(n+2)}\LC(1\!-\!a)^{-\F{n}{2}}
\LR a^2n^2+2a^2n-an^2-14an-24a+48\RR\RD\\
& &
\LD-an^3-an^2-10an+24a-48\RC,\\
C_5&=&\F{1}{6a(n-2)n}\LC-(1\!-\!a)^{-\F{n}{2}}\LR a^2n^2+4a^2n-an^2-10an-36a+48\RR\RD\\
& &\LD+an^3-3an^2+14an-36a+48\RC.
\ean
$C_i$ for the trace in the dimension $n=2$:
\ban
C_1&=&\F{2-a}{2(1\!-\!a)},\quad C_2\ =\ \F{2+a}{6(1\!-\!a)},\quad C_3\ =\ \F{1}{12},
\quad C_4\ =\ -\F{1}{6},\\
C_5&=&-\F{(a-4)\ln(1\!-\!a)}{2a}-\F{11a-14}{6(1\!-\!a)}.
\ean

\section{Conclusion}
The program computes $E_2$ with torsion for operator
(\ref{nmini}) rather easily (about 10 sec  on a Pentium-75 PC).
Unfortunately, computational complexity of the problem under consideration
is very high. For example, the timings for torsionless computations of
$E_2$ and $E_4$ for the same operator (\ref{nmini}) are $< 1$ sec
and 4 h 5 min, correspondingly. It is clear that the inclusion of the
torsion increases the computational efforts considerably and the
computation of $E_4$ with torsion may take too much time.
Another problem is the volume of
the resulting expressions. There are two ways to handle this
problem. First of all, some work is needed for developing of
algorithms for further reduction of large tensor expressions.
However, due to the natural complexity of the heat invariants,
one can not hope to make the higher-order invariants tractable
by hand. Thus, the methods for automatic usage of these invariants
should be elaborated.

\section*{Acknowledgements}
I would like to thank V. Gusynin for initiating this work
and helpful communications.
This work was supported in part by INTAS project No. 96-0842 and
RFBR project No. 98-01-00101.


\end{document}